\documentclass{article}

\usepackage{microtype}
\usepackage{graphicx}
\usepackage{subcaption}
\usepackage{booktabs}

\usepackage{hyperref}

\usepackage[preprint]{icml2026}

\usepackage{amsmath}
\usepackage{amssymb}
\usepackage{amsfonts}
\usepackage{mathtools}
\usepackage{amsthm}

\usepackage{algorithm}
\usepackage{algorithmic}
\usepackage{float}
\usepackage{tabularx}
\usepackage{array}
\usepackage{nicefrac}
\usepackage{url}

\usepackage[capitalize,noabbrev]{cleveref}

\theoremstyle{plain}

\theoremstyle{definition}

\theoremstyle{remark}

\usepackage[textsize=tiny]{todonotes}

\icmltitlerunning{OptiLoop: Coordination-in-the-Loop Verification and Repair for LLM-Generated Optimization Agents}

\begin{document}

\twocolumn[
  \icmltitle{OptiLoop: Coordination-in-the-Loop Verification and Repair \\
  for LLM-Generated Optimization Agents}

  \begin{icmlauthorlist}
    \icmlauthor{Yujia Xu}{amazon}
    \icmlauthor{Zhiheng Wang}{amazon}
    \icmlauthor{Thi Dinh}{amazon}
  \end{icmlauthorlist}

  \icmlaffiliation{amazon}{Amazon, Seattle, WA 98109, USA}

  \icmlcorrespondingauthor{Yujia Xu}{yujiaxu@amazon.com}
  \icmlcorrespondingauthor{Zhiheng Wang}{zhihengw@amazon.com}
  \icmlcorrespondingauthor{Thi Dinh}{thid@amazon.com}

  \icmlkeywords{Large Language Models, Optimization, Verification, Decentralized Coordination, ADMM}

  \vskip 0.3in
]

\printAffiliationsAndNotice{}  

\begin{abstract}
Many decentralized decision problems require multiple parties to coordinate on shared decisions while keeping objectives, constraints, and data private. Large language models (LLMs) offer a promising way to lower the barrier to participation by generating local optimization agents from natural-language specifications. In coordination settings, however, executability is not enough: a generated agent may compile, solve, and pass local checks while still being semantically wrong, for example by misrepresenting costs, mis-scoping constraints, or responding incorrectly to incentives. Such errors often surface only during coordination, as systematic behavioral failures rather than infeasibility. We propose coordination-in-the-loop verification and repair for LLM-generated optimization agents. We instantiate this idea with an Alternating Direction Method of Multipliers (ADMM)-style consensus protocol and introduce OptiLoop, a pipeline that generates local optimization agents from text, verifies them through short, bounded coordination runs against a fixed reference counterparty, extracts structured behavioral and static evidence, and applies evidence-driven repair. When failures are structural rather than implementational, OptiLoop escalates from localized code fixes to corrected-formulation repair, and it can additionally reuse episodic lessons from prior instances. On 40 held-out test scenarios, OptiLoop-Full improves objective match from 66.0\% to 93.0\% and social match from 68.5\% to 89.0\% relative to a strong local-validation baseline, while reducing mean objective gap from 15.3\% to 3.5\% and mean social gap from 7.6\% to 2.0\%. These results show that, for generated optimization agents deployed inside decentralized decision loops, correctness should be validated in the loop itself rather than through isolated execution alone.
\end{abstract}

\section{Introduction}
\begin{figure*}[t]
  \centering
  \includegraphics[width=\textwidth]{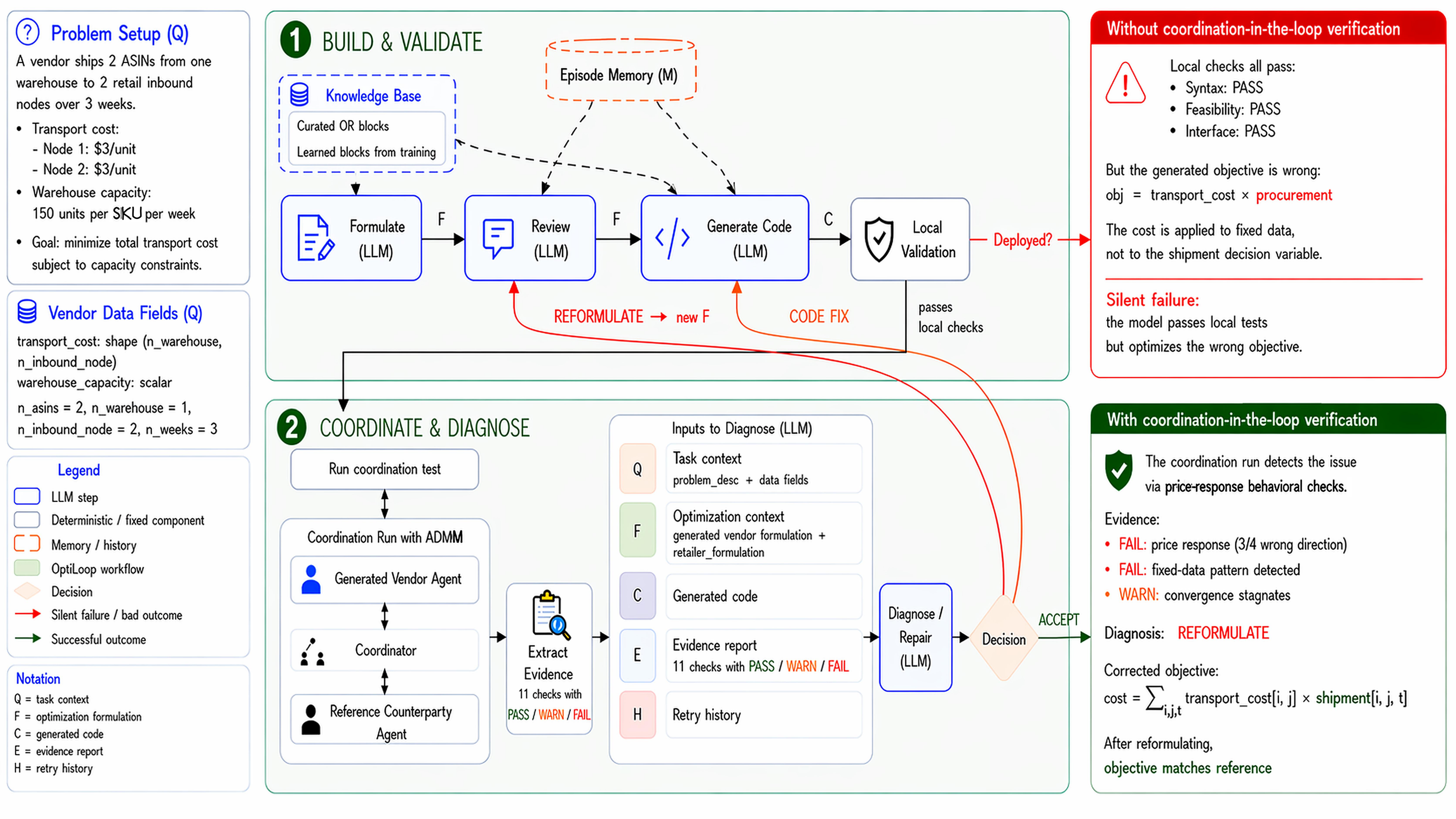}
  \caption{OptiLoop coordination-in-the-loop verification and repair workflow.}
  \label{fig:workflow}
\end{figure*}

Many decision problems of practical interest are decentralized: multiple stakeholders hold private objectives, constraints, and data, yet must still coordinate on shared decisions. This setting arises in supply planning, resource allocation, and scheduling, where privacy requirements, organizational boundaries, and stakeholder autonomy make a single monolithic optimization model undesirable or infeasible. A common alternative is to use iterative coordination methods, such as Alternating Direction Method of Multipliers (ADMM)-style consensus protocols, which repeatedly solve local subproblems while updating shared decisions through structured feedback signals such as dual prices, penalty terms, and residual-based diagnostics \cite{nedic2009distributed,boyd2011admm}. These protocols can scale well while respecting decentralized control, but they rely on a central assumption: each participant implements a correct local optimization agent that faithfully represents its intended private objective and constraints.

In practice, this assumption is often the primary bottleneck. Implementing a local optimization agent requires operations-research expertise to formalize business logic as variables, constraints, and objectives, and engineering expertise to implement solver code that correctly responds to protocol feedback. Large language models (LLMs) offer a potential way to reduce this burden by translating natural-language specifications into mathematical formulations and executable optimization code. Prior work has shown that LLMs can often generate executable optimization models from text \cite{ramamonjison2023nl4opt,ahmaditeshnizi2024optimus,huang2025orlm}. In decentralized coordination settings, however, executability alone is insufficient.

A generated optimization agent may compile, solve, and pass local checks while still being semantically wrong. It may omit a cost term, attach a variable cost to fixed data, or mis-scope a constraint while still returning feasible, well-formed outputs. Such failures often do not manifest as infeasibility; instead, they alter how the agent responds to coordination incentives. As a result, errors that remain invisible under isolated execution may appear only inside the coordination loop, as wrong-sign price response, degenerate proposals, slow or unstable convergence, or persistent disagreement with a trusted counterparty. Thus, isolated execution can verify that an agent runs, but not that it behaves coherently under coordination.

This paper proposes coordination-in-the-loop verification and repair for LLM-generated optimization agents. The central idea is to validate a candidate agent using the structured feedback signals produced by the coordination protocol itself, treating coordination trajectories as semantics-aware verification evidence. We instantiate this principle with an ADMM-style consensus protocol and introduce \textbf{OptiLoop}, a pipeline for generating, verifying, and repairing optimization agents from natural language. As shown in Figure~\ref{fig:workflow}, OptiLoop first generates and locally validates a candidate agent, then runs a short coordination episode against a fixed trusted counterparty to extract behavioral and static evidence, and finally uses that evidence to select among acceptance, localized code repair, and corrected-formulation repair. The key difference from execution-driven validation is that coordination is used not merely to confirm that an agent runs, but to determine whether it behaves correctly under incentives and to drive the choice of repair. We study this problem in a stack-relative setting, where the coordinator and trusted counterparty are fixed deployment infrastructure and the generated agent must be verified for correct behavior within that specific environment.

\paragraph{Contributions.}
First, we formalize coordination-in-the-loop verification as a validation principle for LLM-generated local optimization agents in decentralized coordination settings, emphasizing deployment-time feedback signals rather than isolated execution checks. Second, we introduce OptiLoop, which operationalizes this principle in an ADMM-style consensus protocol via diagnostic coordination episodes, structured evidence extraction, and escalation from localized code repair to corrected-formulation repair for structural failures. Third, we provide a benchmark and ablation study showing that coordination-in-the-loop verification exposes semantic failures missed by local validation and materially improves correctness on held-out decentralized scenarios.

\section{Related Work}
\label{sec:related_work}

OptiLoop sits at the intersection of three lines of work: iterative coordination for decentralized optimization, LLM-based optimization modeling from natural language, and verification and repair for generated programs. Our contribution is to connect these strands through coordination-in-the-loop verification: instead of validating a generated optimization agent only in isolation, we evaluate and repair it using the structured feedback signals produced by the coordination protocol itself.

\paragraph{Iterative coordination and decomposition.}
Large-scale and decentralized optimization problems are often solved by decomposing them into local subproblems and coordinating shared decisions through structured feedback. Distributed optimization methods use message passing to align local objectives \cite{nedic2009distributed}. Consensus and penalty-based methods such as ADMM coordinate local copies through dual variables and proximal penalties, with convergence monitored through residuals \cite{boyd2011admm}. Dual and Lagrangian decomposition similarly update multiplier-like prices in response to constraint violations \cite{fisher1981lagrangian}, while master--subproblem methods such as Benders decomposition refine a master problem using cuts returned by subproblems \cite{rahmaniani2017benders, xu2026dynamic}. Across these methods, a common pattern is iterative interaction through feedback signals. OptiLoop builds on this structure by treating coordination trajectories as a source of verification evidence.

\paragraph{LLMs for optimization modeling and code generation.}
A growing literature studies how to translate natural-language problem descriptions into optimization models and implementations. NL4Opt frames this problem as a benchmark for extracting variables, constraints, and objectives from text \cite{ramamonjison2023nl4opt}. OptiMUS proposes a modular workflow that formulates optimization models, generates solver code, executes it, and iteratively debugs failures \cite{ahmaditeshnizi2024optimus}. ORLM studies automated optimization modeling using synthetic instruction data and benchmark suites \cite{huang2025orlm}. Related work also explores multi-agent LLM reasoning for operations research tasks \cite{xiao2023chain, quan2025leveraging} and structure-guided modeling procedures \cite{pan2025guiding}. Together, these systems show that LLMs can often generate executable optimization models, especially when paired with modular prompting, decomposition, or iterative refinement. 

\paragraph{Verification and repair.}
Most LLM-based optimization-modeling systems still validate candidates primarily through local signals such as successful execution, feasibility, and basic correctness checks, then rely on debugging or regeneration when failures are detected \cite{ahmaditeshnizi2024optimus}. A complementary line of work studies feedback-driven self-correction with episodic memory, showing that language agents can improve across attempts by storing and reusing lessons without updating model weights \cite{shinn2023reflexion}. Recent work has also emphasized that local validation is incomplete for optimization because code may execute successfully while the underlying formulation remains semantically wrong \cite{lian2026reloop}. Other approaches introduce optimization-specific self-correction objectives \cite{jiang2024llmopt} or reduce long repair loops through batched generation and selection \cite{bouscary2026optihive}.

OptiLoop is closest in spirit to recent efforts on optimization-aware verification and repair, but differs in the source of its verification signal. Rather than relying only on isolated execution, solver-side checks, or generic self-correction, OptiLoop uses short, bounded coordination runs to extract structured evidence from incentive response, agreement dynamics, and convergence behavior. It then escalates from localized code fixes to corrected-formulation repair when failures are structural.

\section{Problem Setting}
\label{sec:problem_setting}

We consider a decentralized setting in which multiple parties must coordinate on a shared public decision while keeping their objectives, constraints, and auxiliary variables private. Let \(x \in \mathbb{R}_{\ge 0}^{A \times J \times T}\) denote the public decision tensor, indexed by item \(i \in \{1,\dots,A\}\), node \(j \in \{1,\dots,J\}\), and period \(t \in \{1,\dots,T\}\). Party \(m\) has private data \(\Theta_m\), private auxiliary variables \(y_m\), and a local copy \(x_m\) of the public decision. Its local optimization problem is
\begin{equation}
\min_{x_m,\,y_m}\; c_m(x_m,y_m;\Theta_m)
\quad \text{s.t.} \quad (x_m,y_m)\in \mathcal{F}_m(\Theta_m).
\label{eq:local_model}
\end{equation}
Define the induced private value function
\begin{equation}
f_m(x) := \min_{y_m}\; c_m(x,y_m;\Theta_m)
\quad \text{s.t.} \quad (x,y_m)\in \mathcal{F}_m(\Theta_m),
\label{eq:value_function}
\end{equation}
and the party-specific feasible public-decision set
\[
\begin{aligned}
\mathcal{X}_m := \{x \in \mathbb{R}_{\ge 0}^{A \times J \times T} \;:\;& \exists\, y_m \text{ such that}\\
& (x,y_m)\in \mathcal{F}_m(\Theta_m)\}.
\end{aligned}
\]
For \(M\) parties, a canonical consensus problem is
\begin{equation}
\begin{aligned}
\min_{x_1,\dots,x_M}\;& \sum_{m=1}^M f_m(x_m) \\
\text{s.t.}\;& x_1=\cdots=x_M, \quad x_m \in \mathcal{X}_m \quad \forall m.
\end{aligned}
\label{eq:consensus}
\end{equation}
We instantiate coordination with an ADMM-style consensus protocol \cite{boyd2011admm}. The coordinator maintains a consensus plan \(z \in \mathbb{R}_{\ge 0}^{A \times J \times T}\), party-specific dual tensors \(\lambda_m \in \mathbb{R}^{A \times J \times T}\), and elementwise penalty weights \(\rho \in \mathbb{R}_{>0}^{A \times J \times T}\). At iteration \(k\), party \(m\) solves
\begin{equation}
\begin{aligned}
x_m^{k+1} \in & \arg\min_{x \in \mathcal{X}_m}\;  f_m(x) \\
& - \langle \lambda_m^k, x \rangle 
+ \tfrac{1}{2}\|\sqrt{\rho^k}\odot(x-z^k)\|_2^2.
\end{aligned}
\label{eq:admm_local}
\end{equation}
where \(\odot\) denotes elementwise multiplication. As in standard consensus ADMM, the linear dual term encodes the current coordination incentives, while the quadratic penalty term stabilizes the local update by penalizing deviation from the current consensus iterate. The coordinator then updates
\begin{align}
z^{k+1} &= \frac{1}{M}\sum_{m=1}^M x_m^{k+1},
\label{eq:admm_z_update}
\\
\lambda_m^{k+1} &= \lambda_m^k + \rho^k \odot \bigl(z^{k+1} - x_m^{k+1}\bigr),
\qquad m=1,\dots,M.
\label{eq:admm_lambda_update}
\end{align}
We monitor the primal and dual residuals
\begin{align}
r^{k+1} &= \Big(\sum_{m=1}^M \|x_m^{k+1}-z^{k+1}\|_2^2\Big)^{1/2},
\label{eq:primal_residual}
\\
s^{k+1} &= \sqrt{M}\,\|\rho^k \odot (z^{k+1}-z^k)\|_2,
\label{eq:dual_residual}
\end{align}
and terminate when both fall below prescribed tolerances. These residuals measure, respectively, disagreement among local copies of the public decision and progress of the consensus iterate across coordination steps. Exact stopping criteria and penalty-adaptation rules are deferred to Appendix~\ref{app:admm_details}.

OptiLoop treats each party as a black-box optimization agent implementing
\[
\begin{aligned}
\texttt{solve}(z^k,\lambda_m^k,\rho^k)
&\;\longrightarrow\;
\bigl(x_m^{k+1},\; \texttt{obj\_decomp}_m^{k+1}\bigr).
\end{aligned}
\]
Here \(x_m^{k+1}\) is the proposed public plan and \(\texttt{obj\_decomp}\) is an objective decomposition aligned with \eqref{eq:admm_local}. In our implementation,
\begin{equation}
\begin{aligned}
&\texttt{private}(x) \approx f_m(x),\\
&\texttt{price}(x) = \langle \lambda_m, x\rangle,\\
&\texttt{prox}(x) = \tfrac{1}{2}
  \bigl\|\sqrt{\rho}\odot(x-z)\bigr\|_2^2.
\end{aligned}
\label{eq:obj_decomp}
\end{equation}
so that the local augmented objective can be read as
\[
\texttt{private}(x) - \texttt{price}(x) + \texttt{prox}(x).
\]
Here, \(\texttt{price}(x)\) captures the current coordination incentive induced by the dual signal: increasing a component of \(x\) is rewarded or penalized according to the corresponding entry of \(\lambda_m\). The term \(\texttt{prox}(x)\) is a quadratic penalty that discourages large deviations from the current consensus plan \(z\), with \(\rho\) controlling the strength of this pull component-wise.

In our target application, correctness is defined relative to a fixed deployed coordination stack: the coordinator, interface contract, and trusted counterparty implementation are held fixed, while the generated agent is evaluated under the scenario-specific inputs and coordination signals that arise within that environment. Under this definition, correctness is not an abstract property of the generated code in isolation; it is a behavioral property of the agent within the deployed coordination environment.

We do not propose a new coordination algorithm; the ADMM-style loop and the reference counterparty are fixed infrastructure, and our contribution is to generate, verify, and repair a local optimization agent so that it behaves correctly under feedback-driven coordination.

\section{Methodology}
\label{sec:methodology}

\subsection{OptiLoop Pipeline}
\label{subsec:overview}

Given a natural-language problem description and typed data fields, OptiLoop generates an executable local optimization agent implementing the interface defined in Section~\ref{sec:problem_setting}. The pipeline has four stages: (1) formulation and code generation, (2) lightweight local validation, (3) in-loop verification through a short, bounded coordination run against a fixed counterparty, and (4) evidence-driven diagnosis and repair. The key idea is that coordination behavior provides a semantics-aware verification signal: an agent that appears plausible under isolated execution may still behave incorrectly when exposed to coordination incentives. This bounded coordination run is not used merely to confirm executability, but to test whether the generated agent responds coherently to prices, penalties, and consensus pressure inside the deployed feedback loop. Formulation and code generation are knowledge-augmented through a curated library of operations-research building blocks, together with learned extensions distilled from training scenarios, that provide reusable coordination, constraint, and objective templates; additional details are deferred to Appendix~\ref{app:knowledge_base}.

A full pseudocode listing of the pipeline is given in Algorithm~\ref{alg:optiloop_appendix}.
Here, \(q\) denotes the natural-language problem specification, \(\mathcal{D}\) the typed data schema, \(\mathcal{A}_{\mathrm{fix}}\) the fixed trusted counterparty, and \(R\) the maximum repair budget. OptiLoop maintains a structured formulation \(F\), generated solver code \(C\), optional retrieved memory \(\mathcal{M}\), a coordination trajectory \(\tau\), an extracted evidence report \(E\), and a diagnosis action \(a\). The algorithm shows OptiLoop's bounded repair loop: after a failed diagnosis, it applies either localized code repair or corrected-formulation repair, up to budget \(R\).
\begin{algorithm*}[t]
\caption{OptiLoop: Coordination-in-the-Loop Verification and Repair}
\label{alg:optiloop_appendix}
\centering
\begin{minipage}{0.9\textwidth}
\begin{algorithmic}[1]
\REQUIRE Natural-language specification $q$, typed data schema $\mathcal{D}$, fixed counterparty agent $\mathcal{A}_{\mathrm{fix}}$, maximum repair budget $R$
\STATE $\mathcal{M} \leftarrow \textsc{RetrieveMemory}(q)$ \COMMENT{optional; empty if memory is disabled}
\STATE $F \leftarrow \textsc{Formulate}(q, \mathcal{D}, \mathcal{M})$ \COMMENT{structured mathematical formulation}
\STATE $C \leftarrow \textsc{CodeGen}(F, \mathcal{D}, \mathcal{M})$ \COMMENT{solver code implementing \texttt{solve}}

\FOR{$r = 1,2,\ldots,R$}
    \STATE $\textsf{local} \leftarrow \textsc{LocalValidate}(C)$
    \IF{$\textsf{local} = \textsc{FailHard}$}
        \STATE $C \leftarrow \textsc{PatchOrRegenerate}(C, \textsf{local})$
        \STATE \textbf{continue}
    \ENDIF

    \STATE $\tau \leftarrow \textsc{CoordRun}(C, \mathcal{A}_{\mathrm{fix}})$ \COMMENT{short coordination episode}
    \STATE $E \leftarrow \textsc{ExtractEvidence}(\tau, C)$ \COMMENT{behavioral and static checks}
    \STATE $a \leftarrow \textsc{Diagnose}(q, F, C, E)$ \COMMENT{$a \in \{\textsc{Accept}, \textsc{CodeFix}, \textsc{Reformulate}\}$}

    \IF{$a = \textsc{Accept}$}
        \STATE \textbf{return} $(F, C)$
    \ELSIF{$a = \textsc{CodeFix}$}
        \STATE $C \leftarrow \textsc{CodeFix}(C, E)$ \COMMENT{localized implementation repair}
    \ELSIF{$a = \textsc{Reformulate}$}
        \STATE $F \leftarrow \textsc{CorrectFormulation}(q, F, E, \mathcal{M})$ \COMMENT{structural formulation repair}
        \STATE $C \leftarrow \textsc{CodeGen}(F, \mathcal{D}, \mathcal{M})$
    \ENDIF
\ENDFOR

\STATE \textbf{return} $(F, C)$ \COMMENT{best candidate after bounded repair}
\end{algorithmic}
\end{minipage}
\end{algorithm*}

\subsection{Evidence Extraction from Coordination}
\label{subsec:evidence}

Local validation alone cannot reliably detect semantic modeling errors that preserve executability while distorting economic behavior. OptiLoop addresses this gap by converting a short coordination episode into a structured evidence report \(E\) for diagnosis and repair. We run the generated agent against a fixed counterparty, observe the trajectory \(\tau\) of proposals, consensus variables, dual signals, penalties, residuals, and objective decompositions, and combine these in-loop signals with lightweight static analysis of the generated code.

The evidence report contains two complementary sources of evidence. Behavioral checks evaluate whether the agent responds coherently to coordination signals, including incentive response, degeneracy, marginal behavior, and agreement dynamics. Static checks target code- and formulation-level anti-patterns, such as missing cost terms, costs applied to fixed data rather than decision variables, inconsistent private decision variables, or malformed outputs. Together, these checks address the common case in which a generated model passes local validation but encodes the wrong economics or constraint structure. Static checks can localize suspicious implementation patterns, while in-loop behavioral checks test whether the agent behaves correctly as a participant in the coordination protocol.

Table~\ref{tab:checks_list} summarizes the evidence checks used by OptiLoop, listing each check's source, severity, and diagnostic signal; implementation thresholds are deferred to Appendix~\ref{app:checks}.

\begin{table*}[tb]
\centering
\small
\setlength{\tabcolsep}{4pt}
\renewcommand{\arraystretch}{1.08}
\begin{tabularx}{0.98\textwidth}{@{}c
  >{\raggedright\arraybackslash}p{3.3cm}
  c
  c
  >{\raggedright\arraybackslash}X@{}}
\toprule
\# & \textbf{Check} & \textbf{Source} & \textbf{Severity} & \makebox[\linewidth][c]{\textbf{Diagnostic signal}} \\
\midrule
1  & Utility sign consistency
   & Behavioral & Hard
   & Wrong sign in the private utility convention. \\

2  & Price-response sanity test
   & Behavioral & Hard
   & PO moves in the wrong direction under positive coordination prices. \\

3  & Convergence trajectory
   & Behavioral & Soft
   & Residuals converge, decrease, stagnate, oscillate, or decay slowly. \\

4  & Degenerate plan
   & Behavioral & Hard / warning
   & All-zero or near-zero PO plan when nontrivial fulfillment is expected. \\

5  & Social-gradient perturbation
   & Behavioral & Informational
   & Local mismatch between private response and social objective near consensus. \\

6  & Cost--data anti-pattern
   & Static & Hard / warning
   & Cost parameter is applied to fixed data instead of a decision variable. \\

7  & Decision-variable audit
   & Static & Informational
   & Private decision variables are checked against the described decisions. \\

8  & Marginal-cost consistency
   & Behavioral & Soft
   & Marginal utility changes in the wrong direction as PO increases. \\

9  & Demand-coverage sanity
   & Behavioral & Informational
   & Proposed decisions appear inconsistent with demand or service requirements. \\

10 & Dual-price outlier sanity
   & Behavioral & Informational
   & Large or unstable dual prices indicate persistent disagreement. \\

11 & Missing cost parameters
   & Static & Hard
   & Cost-like input field appears in data but is absent from generated code. \\
\bottomrule
\end{tabularx}
\caption{
Evidence checks used by OptiLoop. Behavioral checks use solver outputs or coordination traces, while static checks analyze the generated formulation or code without running coordination. Hard signals require repair or escalation; soft signals contribute supporting evidence; informational signals help localize the failure but do not by themselves trigger repair. The evidence-class ablation in Appendix~\ref{app:evidence_ablation} evaluates static and behavioral evidence independently.
}
\label{tab:checks_list}
\end{table*}

\subsection{Diagnosis and Repair}
\label{subsec:diagnosis_repair}
Given the evidence report \(E\), OptiLoop selects one of three actions: \textsc{Accept}, \textsc{CodeFix}, or \textsc{Reformulate}. If the agent appears behaviorally coherent and no severe issues are detected, it is accepted. If the evidence indicates a localized implementation problem, OptiLoop applies a targeted code patch. If the evidence instead suggests a structural modeling error, OptiLoop escalates to reformulation. To ensure safety and consistency, the diagnosis output is filtered through simple rules; for example, \textsc{Accept} is disallowed when hard-fail evidence is present, in which case the system escalates to \textsc{Reformulate}. 

\textsc{CodeFix} targets localized implementation failures such as missing fields, shape mismatches, or straightforward parameter-binding errors. Patches are intentionally minimal in order to reduce the chance of introducing new errors. \textsc{Reformulate} addresses cases in which the generated formulation itself is wrong. OptiLoop produces a revised formulation \(F'\) conditioned on the original problem description, the current formulation, and the evidence report, and then regenerates solver code from \(F'\). This allows the system to overwrite faulty modeling choices directly rather than relying on prompt-level retrying alone. Additional implementation details for the diagnosis actions, deterministic post-processing rules, and corrected-formulation repair procedure are provided in Appendix~\ref{app:diagnosis}.

\subsection{Episodic Experience Reuse}
\label{subsec:memory}

OptiLoop optionally uses episodic experience reuse to reduce recurring error patterns. At generation time, it retrieves relevant lessons from prior training instances and injects them into formulation and code generation. The memory store includes both instance-level bug-and-fix patterns and higher-level modeling principles. Memory is updated only on training scenarios and frozen during evaluation. Appendix~\ref{app:memory} provides further details on what is stored in episodic memory, how lessons are retrieved, and how they are injected during generation.

\section{Experiments}
\label{sec:experiments}

\noindent\textbf{Overview.}
We evaluate OptiLoop in a controlled ablation study designed to isolate the contributions of local validation, coordination-in-the-loop verification and repair, and episodic experience reuse. The ablation ladder is cumulative: \textbf{Baseline-Gen} performs natural-language-to-formulation-to-code generation only; \textbf{Baseline-LocalVal} adds syntax checks and local solve validation; \textbf{OptiLoop-Core} further adds coordination-based evidence extraction, diagnosis, and repair; and \textbf{OptiLoop-Full} additionally incorporates frozen episodic memory learned from the training scenarios. All methods use the same underlying LLM and solver stack, and differ only in which verification, repair, and memory components are enabled. These baselines are controlled component ablations of a shared generation pipeline, designed to isolate the incremental value of each OptiLoop component.

\paragraph{Scenarios and split.}
We evaluate on a bounded hand-authored benchmark of 80 scenarios, each defined by a natural-language description, typed data-field metadata, and a concrete data dictionary. Each scenario specifies the generated agent's private modeling problem, including public coordination variables, auxiliary private variables, parameters, constraints, and scenario-specific business rules. The scenarios instantiate recurring supply-chain primitives such as order aggregation, inventory balance, warehouse and transportation capacity, production and shipping costs, holding costs, backlog penalties, and service-level constraints, combined with linear, quadratic, and piecewise cost structures. Appendix~\ref{app:worked_examples} provides worked examples of representative scenarios and failure modes targeted by OptiLoop.

We use 40 scenarios for training and 40 disjoint scenarios for testing. The split is by scenario instance: both splits share common modeling primitives, but the test set introduces harder compositions, less common cost structures, more ambiguous natural-language phrasing, and more complicated constraints, including shared capacity, returns or rework flows, and cross-ASIN coupling. As summarized in Appendix~\ref{app:reproducibility}, the held-out test set also introduces structural tags absent from training, including novel cost structures and deliberate natural-language challenges. Episodic memory is built only from the training scenarios and frozen at test time.

\paragraph{Counterparty and coordination protocol.}
Our evaluation targets an asymmetric onboarding setting in which a newly generated agent is integrated into an existing deployed coordination stack. All coordination episodes therefore pair the generated agent with a hand-coded counterparty agent representing trusted platform-side infrastructure. The coordinator, counterparty implementation, and interface contract are held constant across scenarios, while the prompt-specified inputs and data vary by instance. This design isolates errors in the generated agent while matching the intended deployment setting. In-loop verification uses the ADMM-style consensus protocol from Section~\ref{sec:problem_setting}: during generation-time verification, OptiLoop runs short, bounded coordination episodes to extract evidence from the resulting agent--counterparty interaction, whereas final evaluation uses longer runs with early stopping and adaptive penalty updates.

\paragraph{Metrics.}
We report four metrics. \textbf{Objective match} (\texttt{obj\_match}) measures whether the generated agent reproduces the intended private objective behavior at fixed consensus inputs: for reference consensus plans, we compare its objective value with that of the hand-coded reference agent and count a scenario as matched if the relative error is below 0.1\%. \textbf{Objective gap} is the corresponding mean relative error, averaged across test conditions and capped at 100\% per scenario. \textbf{Social match} (\texttt{social\_match}) measures whether the final coordinated solution is within 0.1\% of the optimum of the corresponding centralized joint problem, whose objective is the sum of the two parties' private objectives. \textbf{Social gap} is the corresponding mean relative optimality gap, again capped at 100\% per scenario. Match metrics capture strict correctness, while gap metrics capture the severity of remaining errors. Objective metrics evaluate whether the generated local agent encodes the intended private model, while social metrics evaluate whether that local correctness translates into a good coordinated outcome.

\paragraph{Implementation and reproducibility.}
All methods use the same underlying LLM, prompting interface, solver stack, and evaluation harness, so performance differences reflect the enabled verification and repair components rather than changes in generation infrastructure. We report five independent runs per method. Although decoding temperature is fixed at zero, the serving stack can exhibit small run-to-run nondeterminism, so repeated runs measure end-to-end robustness rather than sampling diversity.  Additional implementation and reproducibility details are provided in Appendix~\ref{app:exp_details}.

\section{Results}
\label{sec:results}

\subsection{Main correctness results}

\begin{figure*}[tb]
\centering
\includegraphics[width=\textwidth]{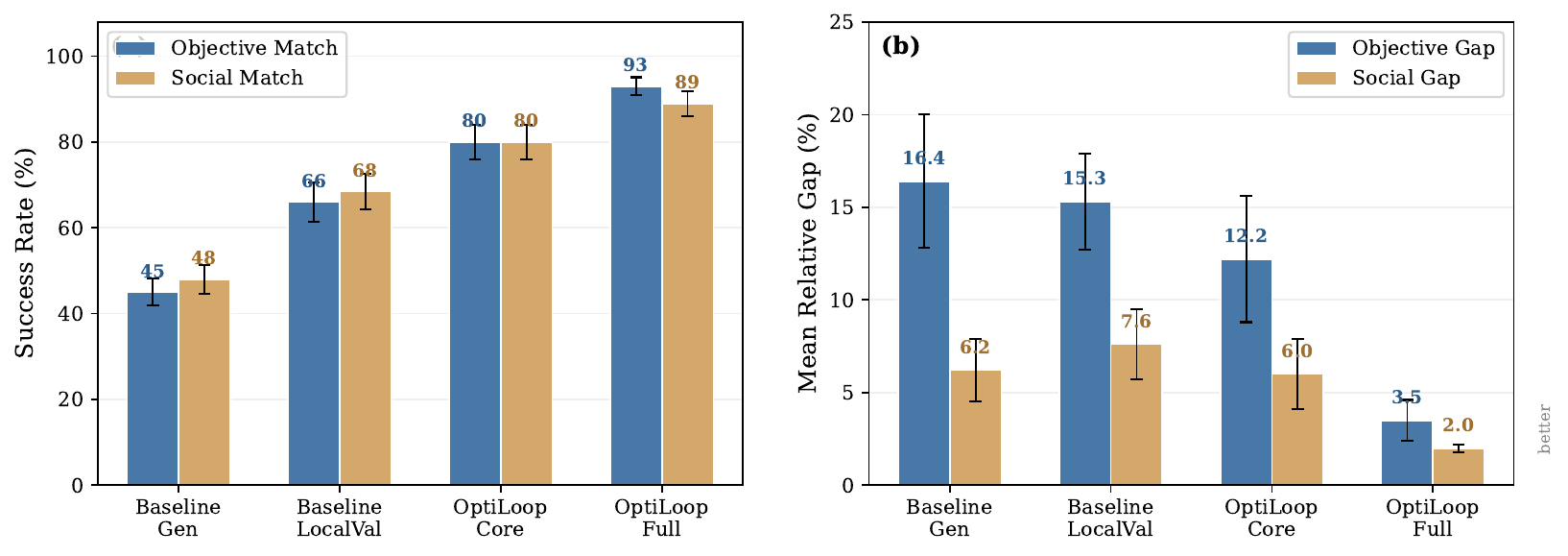}
\caption{Main ablation results on 40 held-out test scenarios. Left: objective match and social match success rates at 0.1\% tolerance. Right: mean relative objective and social gaps to reference solutions, capped at 100\% per scenario.}
\label{fig:results}
\end{figure*}

Figure~\ref{fig:results} summarizes the main correctness results across five independent runs on 40 held-out test scenarios; exact values are reported in Table~\ref{tab:main_results}. The match metrics improve
monotonically across the ablation ladder. Adding local validation to generation improves objective match from 45.0\% to 66.0\% and social match from 48.0\% to 68.5\%, confirming that many failures are still basic execution, feasibility, or interface errors. Adding coordination-in-the-loop diagnosis and repair further improves both match metrics to 80.0\%, indicating that a substantial fraction of locally valid agents still behave incorrectly once exposed to coordination signals.

OptiLoop-Full delivers the strongest overall performance, reaching 93.0\% objective match and 89.0\% social match. The gap metrics show the same pattern, but also quantify the severity of the remaining errors: relative to Baseline-LocalVal, OptiLoop-Full reduces the mean objective gap from 15.3\% to 3.5\% and the mean social gap from 7.6\% to 2.0\%. Thus, coordination-in-the-loop verification improves not only strict pass/fail correctness, but also the quality of imperfect solutions when exact matching is not achieved.

\begin{table*}[tb]
\centering
\small
\begin{tabular}{lcccc}
\toprule
\textbf{Method} & \textbf{Obj Match (\%)} & \textbf{Social Match (\%)} & \textbf{Obj Gap (\%)} & \textbf{Social Gap (\%)} \\
\midrule
Baseline-Gen          & $45.0 \pm 3.1$ & $48.0 \pm 3.3$ & $16.4 \pm 3.6$ & $6.2 \pm 1.7$ \\
Baseline-LocalVal     & $66.0 \pm 4.5$ & $68.5 \pm 4.2$ & $15.3 \pm 2.6$ & $7.6 \pm 1.9$ \\
OptiLoop-Core         & $80.0 \pm 4.0$ & $80.0 \pm 4.0$ & $12.2 \pm 3.4$ & $6.0 \pm 1.9$ \\
OptiLoop-Full         & $93.0 \pm 2.1$ & $89.0 \pm 2.9$ & $\;\;3.5 \pm 1.1$ & $2.0 \pm 0.2$ \\
\bottomrule
\end{tabular}
\caption{Correctness metrics (mean $\pm$ std across 5 runs on 40 test scenarios). Match metrics report success rate (\%) at 0.1\% tolerance; gap metrics report mean relative gap to reference (\%), capped at 100\% per scenario.}
\label{tab:main_results}
\end{table*}

\subsection{Ablation insights} \label{subsec:evidence_ablation}

The main ablation results in Figure~\ref{fig:results} separate three sources of improvement. First, Baseline-LocalVal improves over generation alone because many invalid agents fail basic execution, feasibility, or interface checks. Second, OptiLoop-Core improves over Baseline-LocalVal because some agents pass local validation but still behave incorrectly during coordination. These failures are semantic rather than syntactic: the model may solve successfully, but attach costs to fixed inputs, omit relevant cost terms, or respond in the wrong direction to prices. In these cases, the agent is not wrong because it crashes or violates an interface contract; it is wrong because it encodes the wrong economics under coordination incentives. Coordination-in-the-loop checks expose these errors by testing the generated agent as an economic participant rather than only as executable code. Third, OptiLoop-Full further improves over OptiLoop-Core because episodic memory helps prevent recurring mistakes before repair is needed, improving both correctness and residual gaps.

\begin{figure}[tb]
\centering
\includegraphics[width=0.85\columnwidth]{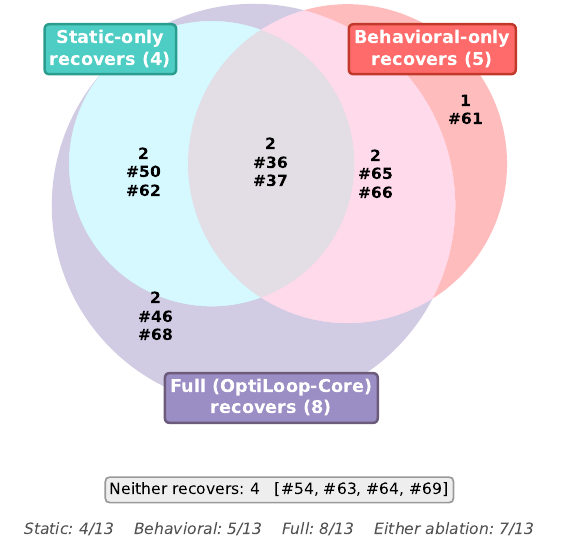}
\caption{ Recovery of the 13 test scenarios where Baseline-LocalVal fails. Each circle represents the set of scenarios recovered by one evidence mode. Numbers inside each region report the number of scenarios in that region, and labels of the form \#k denote benchmark scenario IDs. Static-only and Behavioral-only recover different subsets, and Full recovers additional scenarios that neither ablation fixes alone.
}
\label{fig:evidence}
\end{figure}

We also ablate the evidence used by OptiLoop-Core. Static-only evidence uses code-level and interface checks, Behavioral-only evidence uses coordination traces, and Full combines both. Figure~\ref{fig:evidence} shows the recovery pattern on the 13 test scenarios where Baseline-LocalVal fails to achieve objective match. Static-only recovers 4 scenarios and Behavioral-only recovers 5, but the recovered subsets differ; Full recovers 8 scenarios, including cases that neither ablation fixes alone. This supports the claim that static localization and behavioral confirmation provide complementary signals. Appendix~\ref{app:evidence_ablation} gives the full ablation.

\subsection{Efficiency and repair behavior}

Table~\ref{tab:efficiency} summarizes the efficiency and repair behavior of the pipeline. OptiLoop-Core improves correctness at the cost of additional coordination and diagnosis calls: compared with Baseline-LocalVal, its average cost increases from $3.6$ to $6.7$ LLM calls and from $51$ to $95$ seconds per scenario. This overhead reflects the extra short coordination run, evidence extraction, diagnosis, and possible repair.

OptiLoop-Full recovers much of that overhead through episodic experience reuse. With memory enabled, the average number of LLM calls drops from $6.7$ to $4.9$, runtime decreases from $95$ to $67$ seconds, and the mean number of repair iterations falls from $0.77$ to $0.31$. The diagnosis distribution also shifts: the fraction of scenarios accepted without repair rises from 51.5\% to 80.5\%, while CodeFix becomes rare. This suggests that retrieved lessons primarily prevent recurring mistakes during generation, rather than merely making later repair cheaper.

\begin{table*}[tb]
\centering
\small
\begin{tabular}{lcccccc}
\toprule
\textbf{Method} & \textbf{LLM Calls} & \textbf{Time (s)} & \textbf{Accept} & \textbf{CodeFix} & \textbf{Reformulate} & \textbf{\#Repairs} \\
\midrule
Baseline-Gen          & $2.0$ & $36 \pm 0$ & --- & --- & --- & --- \\
Baseline-LocalVal     & $3.6 \pm 0.2$ & $51 \pm 2$ & --- & --- & --- & --- \\
OptiLoop-Core         & $6.7 \pm 0.2$ & $95 \pm 3$ & $51.5\%$ & $25.5\%$ & $23.0\%$ & $0.77$ \\
OptiLoop-Full         & $4.9 \pm 0.2$ & $67 \pm 3$ & $80.5\%$ & $1.0\%$ & $18.5\%$ & $0.31$ \\
\bottomrule
\end{tabular}
\caption{Efficiency and repair metrics (mean $\pm$ std). \textbf{LLM Calls}: mean LLM invocations per scenario. \textbf{Time}: wall-clock seconds per scenario. \textbf{Accept/CodeFix/Reformulate}: fraction of scenarios ending with each diagnosis action. \textbf{\#Repairs}: mean number of repair iterations per scenario.}
\label{tab:efficiency}
\end{table*}

\subsection{Failure analysis}
\label{subsec:failure_analysis}

Appendix~\ref{app:failure_mode} provides a bug-level taxonomy of the residual failures, mapping each failure mode to its root cause, strongest evidence signal, evidence class, and typical repair action. This taxonomy complements the broader scenario-level categories discussed below by showing which errors are primarily exposed by behavioral checks, which are caught by static checks, and which require both. It also clarifies why in-loop behavioral evidence is necessary: many residual errors are structurally valid programs that remain executable and locally plausible, yet still encode the wrong economics or constraint scope under coordination.

Even with the full pipeline, 7\% of scenario-run instances fail, corresponding to 14 failure instances across 6 distinct scenarios. These residual cases are mostly not syntax or interface failures; they are structural or interpretive errors that can remain plausible under local validation. We group them into the following three broad categories:

\paragraph{Novel cost structures.}
Some failures involve cost forms that are underrepresented in the generation and repair context, especially piecewise or tiered costs such as regular-versus-overtime production. In these cases, the generated formulation can respond reasonably near the tested consensus points while still encoding the wrong global cost decomposition. Short coordination runs expose some incentive errors, but may not fully identify mistakes that appear only in untested operating regions.

\paragraph{Ambiguous specifications.}
A second class of failures comes from natural-language descriptions that admit multiple reasonable formalizations. Examples include underspecified capacity sharing, return or rework flows, and colloquial cost descriptions. In-loop evidence can show that a candidate agent behaves inconsistently with the reference stack, but it cannot always determine which business interpretation was intended without additional clarification.

\paragraph{Complex constraint composition.}
The remaining failures tend to combine several otherwise familiar modeling primitives at once. When inventory balance, order aggregation, shared capacities, lane restrictions, and service constraints interact, the generated agent may omit one constraint, apply it at the wrong index level, or aggregate the wrong variables. Such errors often preserve feasibility and executability, but distort the coordinated solution. These cases suggest that future versions should probe a broader set of operating points during verification and explicitly request clarification when the specification is intrinsically ambiguous.

Appendix~\ref{app:worked_examples} complements the aggregate results with three concrete scenarios. These examples show how the same pipeline behaves across different difficulty levels: a standard case where static checks catch an objective-construction mistake, a coordination-sensitive case where the error is visible only through behavioral response, and a hard case where novel costs or ambiguous language leave multiple plausible formulations. The second case most directly illustrates the motivation for in-loop verification, since the agent appears locally valid but reveals its semantic error only when placed in the coordination protocol.

\section{Discussion and Limitations}
\label{sec:discussion}

OptiLoop is instantiated with an ADMM-style consensus protocol, but its broader contribution is the coordination-in-the-loop verification principle: when generated optimization agents operate inside feedback-driven decision loops, those signals provide semantics-aware evidence for verification and repair. Our current evidence library is therefore specialized to ADMM-style coordination, relying most strongly on price response, residual behavior, and short coordination trajectories. Extending this library to other coordination and decomposition mechanisms, such as Benders decomposition, column generation, or dual decomposition, is an important direction for future work.

Our evaluation also has several limitations. First, it is matched to the target onboarding setting, where a generated agent is verified against a fixed coordinator and trusted counterparty, so transfer across coordination stacks or counterparties remains open. Second, our experiments use a bounded synthetic scenario suite with reference objectives computed by hand-coded solvers, so broader scenario coverage, larger problem scales, and public benchmark suites would strengthen the empirical case for generality. Third, although OptiLoop-Full is substantially more efficient than OptiLoop-Core, it remains more expensive than simple generation or local validation, making the current system better suited to offline onboarding and pre-deployment verification than to latency-sensitive settings.

\section{Conclusion}

We introduced OptiLoop, a coordination-in-the-loop verification and repair pipeline for LLM-generated optimization agents. By using short coordination episodes to extract structured behavioral evidence, OptiLoop supports targeted repair ranging from localized code fixes to corrected-formulation repair for structural modeling errors. In a controlled benchmark study, coordination-in-the-loop diagnosis and episodic experience reuse consistently improve correctness and reduce error relative to local validation alone.

Taken together, these findings show that local validation is incomplete for decentralized LLM-generated optimization agents: some semantic failures remain hidden under isolated execution and surface only during coordination. This suggests that generated optimization agents should be evaluated, and when necessary repaired, in the same feedback context in which they will ultimately be deployed. 

\clearpage

\bibliography{mybibliography}
\bibliographystyle{icml2026}

\clearpage
\appendix
\onecolumn

\section{Additional Methodology Details}
\label{app:method_details}

\subsection{Knowledge-Augmented Generation}
\label{app:knowledge_base}

OptiLoop integrates a curated library of operations-research building blocks into the generation pipeline. The library contains reusable modeling patterns, including coordination terms (e.g., dual-price/subsidy terms and proximal penalties), inventory balance constraints, order aggregation constraints, warehouse capacity limits, truck volume limits, and objective-assembly templates. Each block is represented in both mathematical notation and solver-oriented Python/Xpress form. During code generation, the same blocks, together with a one-shot example and the required \texttt{solve} interface specification, guide the model in producing executable Xpress solver code.

Beyond the curated library, OptiLoop maintains a learned knowledge base that grows during training. After each successful generation, an additional LLM call compares the generated code against the existing block library and extracts any new reusable constraint or objective patterns not already covered, together with common-pitfall notes derived from error-fix pairs encountered during retries. These learned blocks are appended to the curated base and made available in future generations.

This knowledge base is distinct from episodic experience reuse. The OR block library provides reusable modeling primitives during generation, whereas episodic memory retrieves instance-level lessons and repair patterns from prior runs.

\subsection{Evidence checks}
\label{app:checks}

OptiLoop computes an evidence report from two sources: short coordination runs and lightweight static analysis of the generated agent code. Some checks are hard-fail signals that trigger mandatory escalation, while others provide soft or informational evidence used by diagnosis.

\paragraph{Thresholds and constants.}
Unless otherwise noted, the degeneracy warning threshold is 80\% zeros; the price-response test samples 4 public variables and fails if a majority move in the wrong direction; ADMM verification runs are capped at 300 iterations and evaluation runs at 2000 iterations; and the stopping tolerance is $10^{-4}$ for both primal and dual residuals.

\subsection{Diagnosis and repair rules}
\label{app:diagnosis}

\paragraph{Inputs.}
The diagnosis module receives: (i) the natural-language description, (ii) the generated formulation, (iii) the evidence report (Table~\ref{tab:checks_list}), (iv) the generated agent code, (v) data-field shape/type metadata, and (vi) short-term history across retries.

\paragraph{Actions.}
The diagnosis module outputs an action $a \in \{\textsc{Accept},\textsc{CodeFix},\textsc{Reformulate}\}$. We apply deterministic post-processing rules to ensure safety and consistency:
\begin{enumerate}
    \item If $a=\textsc{Accept}$ and \texttt{has\_fail} is true, override to \textsc{Reformulate}.
    \item If $a=\textsc{Accept}$ and \texttt{has\_fail} is false, accept as-is.
    \item If $a=\textsc{CodeFix}$, apply a minimal patch (typically 1--2 lines) and re-run local validation and a short verification episode. If the patch does not change the objective behavior, we accept the patched agent as an escape hatch.
    \item If $a=\textsc{Reformulate}$, perform structural repair as described below; retries are bounded.
\end{enumerate}

\paragraph{Corrected-formulation repair.}
When the action is \textsc{Reformulate}, the repair module may emit an explicit corrected mathematical formulation $F'$ (sets/indices, parameters, variables, objective, constraints). If present, OptiLoop replaces the current formulation with $F'$, regenerates solver code from $F'$, and re-enters verification. This mechanism directly overwrites structural modeling errors and mitigates repeated failure modes observed under prompt-level ``retry with feedback.''

\subsection{Episodic memory}
\label{app:memory}

OptiLoop optionally uses episodic experience reuse to reduce recurring semantic errors. Memory stores (i) scenario-level lessons describing error patterns and fixes and (ii) meta-lessons distilled across scenarios. Retrieved lessons are injected as mandatory rules during formulation and code generation.

\paragraph{What is stored.}
Each lesson includes a brief problem summary, the failure pattern (what went wrong), the fix pattern (what corrected it), and a short ``do/don't'' rule (e.g., ``variable costs must multiply decision variables, not fixed inputs''). Meta-lessons are always included; scenario lessons are retrieved by similarity to the current prompt.

\paragraph{Retrieval and usage.}
At generation time, OptiLoop retrieves the top-$K$ most relevant scenario lessons (we use $K{=}3$) and injects them alongside meta-lessons into the formulation and code-generation prompts. Memory is updated only during training and frozen at test time.

\section{Additional Experimental Details}
\label{app:exp_details}

\subsection{Benchmark split and evaluation setup}
\label{app:reproducibility}

\paragraph{Benchmark split details.}
We evaluate on 80 hand-authored scenarios with a 40/40 train/test split by scenario instance rather than by constraint family. Both splits share common modeling primitives such as order aggregation, inventory balance, warehouse capacity, and transportation costs. The held-out test set, however, introduces additional structural tags, including congestion costs, piecewise costs, shared-capacity constraints, returns or rework flows, cross-ASIN coupling, and deliberately ambiguous or terse natural-language phrasing. This split is intended to evaluate generalization under partial structural overlap: familiar building blocks recur across splits, but the test scenarios require new compositions, less common cost structures, and linguistic interpretations not seen in the scenarios used to construct episodic memory.

\paragraph{Optimization primitives.}
The benchmark is built around a retailer--vendor sourcing setting in which a retail platform coordinates purchase-order decisions with vendors. The shared public decision is the purchase-order quantity, while the generated vendor-side agent models how it can feasibly and profitably respond to that proposed sourcing plan. Depending on the scenario, the vendor-side private model may introduce additional variables such as shipment quantities, warehouse inventory, production quantities, truck counts, shortfall variables, or transfer quantities. These variables determine the vendor's feasible response to the public purchase-order decision.

The constraint structure includes several common families. Order-aggregation constraints link the public purchase-order quantity to private shipment or fulfillment decisions. Inventory-balance constraints track stock across time using initial inventory, inbound flow, outbound flow, and demand or service requirements. Capacity constraints limit warehouse throughput, production, storage, transportation lanes, truck volume, or shared resources. Business-rule constraints encode requirements such as minimum order quantities, blackout periods, allowed lanes, service levels, priority rules, or terminal inventory targets.

The objective terms also vary across scenarios. Vendor models may include transportation costs, production costs, holding costs, backlog or shortfall penalties, fixed truck costs, congestion costs, preference weights, and revenue-like utility terms. These terms are intentionally diverse because many semantic formulation errors arise from objective construction rather than constraint syntax. For example, a generated model may incorrectly treat a decision-dependent transportation cost as a fixed constant, use the wrong sign for a utility term, or omit a penalty that appears only implicitly in the natural-language description.

These primitives make the benchmark more than a code-generation task. The generated vendor-side agent must infer the algebraic role of each data field: whether it is an index set, parameter, decision-dependent coefficient, bound, cost, or target. The held-out scenarios therefore test whether the method can assemble familiar sourcing and supply-chain primitives into new formulations, rather than merely reproduce surface-level code patterns from training scenarios.

\paragraph{Tag overlap.}
Each scenario is annotated with structural tags (e.g., \texttt{order\_aggregation}, \texttt{inventory\_balance}, \texttt{congestion\_cost}). Table~\ref{tab:tag_coverage} summarizes the distribution across splits. Both splits share 22 common tags covering standard modeling primitives, while the test set introduces 28 tags absent from training, including novel cost structures, novel constraint types, and deliberate linguistic challenges. This split therefore evaluates generalization to unseen compositions and cost structures rather than simple memorization of recurring training patterns.

\begin{table}[t]
\centering
\small
\begin{tabular}{lccl}
\toprule
\textbf{Tag} & \textbf{Train} & \textbf{Test} & \textbf{Category} \\
\midrule
\multicolumn{4}{l}{\emph{Shared tags (present in both splits)}} \\
\addlinespace
order\_aggregation & 29 & 36 & Logistics \\
inventory\_balance & 20 & 31 & Inventory \\
holding\_cost & 12 & 19 & Objective \\
transport\_cost & 15 & 11 & Objective \\
warehouse\_capacity & 11 & 9 & Capacity \\
profit\_maximization & 10 & 4 & Objective \\
truck\_volume & 6 & 6 & Logistics \\
safety\_stock & 3 & 6 & Inventory \\
blackout\_window & 4 & 5 & Business rule \\
outbound\_rate & 3 & 4 & Capacity \\
margin\_target & 2 & 4 & Business rule \\
inter\_warehouse\_transfer & 3 & 3 & Logistics \\
\addlinespace
\midrule
\multicolumn{4}{l}{\emph{Representative test-only tags (absent from training)}} \\
\addlinespace
novel\_structure & 0 & 5 & Structure \\
ambiguous\_nl & 0 & 4 & NL challenge \\
multi\_echelon & 0 & 2 & Structure \\
tiered\_cost & 0 & 2 & Cost structure \\
large\_scale & 0 & 2 & Scale \\
max\_composition & 0 & 2 & Complexity \\
factory\_capacity & 0 & 2 & Capacity \\
rephrased & 0 & 2 & NL challenge \\
congestion\_cost & 0 & 1 & Cost structure \\
piecewise\_cost & 0 & 1 & Cost structure \\
shared\_capacity & 0 & 1 & Capacity \\
cross\_asin\_coupling & 0 & 1 & Structure \\
returns\_rework & 0 & 1 & Structure \\
lead\_time & 0 & 1 & Structure \\
inventory\_aging & 0 & 1 & Cost structure \\
terse\_phrasing & 0 & 1 & NL challenge \\
verbose\_phrasing & 0 & 1 & NL challenge \\
\bottomrule
\end{tabular}
\caption{Structural tag distribution across train and test splits. We show the 12 most frequent shared tags and 17 representative test-only tags, out of 22 shared and 28 test-only tags overall. Every scenario also carries the \texttt{QP} tag, since all models are convex quadratic programs.}
\label{tab:tag_coverage}
\end{table}

\paragraph{Two-party setting and counterparty details.}
All experiments use the two-party case, consisting of one generated vendor-side agent and one fixed platform-side counterparty. The counterparty is implemented as a hand-coded Xpress optimization model and represents a canonical retail inventory-control agent. It optimizes over the shared purchase-order plan \(x\), balancing inventory holding cost and backlog penalty subject to standard inventory-balance dynamics, with an additional terminal penalty that encourages ending inventory to remain near a safety-stock target. The counterparty is held fixed across methods, scenarios, and coordination episodes, so that differences in performance reflect the quality of the generated agent and the verification pipeline rather than variation in the trusted platform-side model. Consequently, anomalies observed during in-loop verification can be attributed to the generated agent's formulation or behavior rather than to counterparty variability.

\subsection{ADMM stopping and penalty adaptation.}
\label{app:admm_details}
We instantiate coordination-in-the-loop verification with an ADMM-style consensus protocol \cite{boyd2011admm}. For consistency with the consensus formulation in Section~\ref{sec:problem_setting}, we use the standard consensus-ADMM primal residual and its corresponding dual residual, generalized to elementwise penalty weights.

\paragraph{Stopping rule.}
We monitor the primal and dual residuals
\begin{align}
r^{k+1} &= \Big(\sum_{m=1}^{M}\|x_m^{k+1}-z^{k+1}\|_2^2\Big)^{1/2}, \\
s^{k+1} &= \sqrt{M}\,\|\rho^{k}\odot(z^{k+1}-z^{k})\|_2,
\end{align}
and stop early when both \(r^{k+1} < 10^{-4}\) and \(s^{k+1} < 10^{-4}\).

\paragraph{Penalty adaptation.}
We adapt \(\rho\) to balance primal and dual progress:
\begin{equation}
\rho^{k+1}=
\begin{cases}
2\rho^{k} & \text{if } r^{k+1} > 10\cdot s^{k+1},\\
\tfrac{1}{2}\rho^{k} & \text{if } s^{k+1} > 10\cdot r^{k+1},\\
\rho^{k} & \text{otherwise.}
\end{cases}
\end{equation}
In verification runs we cap ADMM at 300 iterations; in evaluation runs we allow up to 2000 iterations, with the same stopping and adaptation rules.

\subsection{Implementation and evaluation protocol}
\label{app:implementation_eval}
\paragraph{Reproducibility details.}
The LLM is Claude Sonnet 4 (\texttt{us.anthropic.claude-sonnet-4-20250514-v1:0}) accessed via Amazon Bedrock. We use temperature $=0.0$ and a maximum output length of $16{,}000$ tokens. All optimization models are solved with FICO Xpress (QP-capable). When episodic memory is enabled, retrieval uses Amazon Titan Text Embeddings v2 (1024 dimensions). The memory store contains 50 frozen entries derived only from the training split: 40 scenario-specific lessons and 10 meta-lessons. All methods use the same model and inference settings; the reported run-to-run variability arises from small serving-side nondeterminism despite temperature-zero decoding.

\paragraph{Runs and evaluation.}
We report five independent repetitions per method. During generation-time verification, OptiLoop runs short coordination episodes capped at 300 iterations. During final evaluation, coordination runs are capped at 2000 iterations with early stopping and adaptive penalty updates, using the stopping tolerances and penalty-adaptation rules reported in Appendix~\ref{app:admm_details}.

\section{Supplementary Results}
\label{app:supp_results}
\subsection{Additional Evidence-Ablation Details}
\label{app:evidence_ablation}

Section~\ref{subsec:evidence_ablation} summarizes the evidence-type ablation. Here we provide the full ablation setup and scenario-level recovery details. We hold generation, local validation, repair budget, and evaluation protocol fixed, and compare four modes: (i) Baseline-LocalVal, which performs no in-loop diagnosis; (ii) Static-only, which uses only code-level and interface evidence (Checks 6, 7, 11); (iii) Behavioral-only, which uses only in-loop coordination evidence (Checks 1--5, 8--10); and (iv) Full (OptiLoop-Core), which combines both evidence sources. In Static-only mode, the ADMM summary is also suppressed, including residual trajectories, per-variable residuals, and dual prices, so that the diagnosis LLM cannot infer behavioral signals from raw coordination traces.

Figure~\ref{fig:evidence} in the main text visualizes recovery on the 13 test scenarios where Baseline-LocalVal fails to achieve objective match. Static-only evidence recovers 4 scenarios, Behavioral-only evidence recovers 5, and Full evidence recovers 8. The recovered subsets differ: scenarios \#50 and \#62 are recovered only by static checks, which detect missing cost parameters and cost-applied-to-fixed-data patterns, while scenario \#61 is recovered only by behavioral checks, which detect wrong-sign marginal response during coordination. Scenarios \#36 and \#37 are recovered by either evidence type independently. Scenarios \#46 and \#68 are recovered only by Full evidence, indicating that static localization and behavioral confirmation can jointly enable repairs that neither signal enables independently. Scenarios \#54, \#63, \#64, and \#69 remain unrecovered by any method and correspond to the hard residual cases discussed in Appendix~\ref{app:failure_mode}.

Because diagnosis and repair are LLM-mediated, Full evidence is not a strict per-scenario union of Static-only and Behavioral-only evidence. Adding evidence can change the selected repair action, the regenerated formulation, and the resulting coordination trajectory.

\subsection{Failure-mode Taxonomy}
\label{app:failure_mode}

Section~\ref{subsec:failure_analysis} summarizes the residual failure categories. Table~\ref{tab:taxonomy} provides a bug-level taxonomy that maps each failure mode to its root cause, strongest evidence signal, evidence class, and typical repair action. The taxonomy clarifies why local validation is insufficient: several failures produce executable and feasible solver code while still encoding the wrong incentive response, cost structure, or constraint scope.

\begin{table*}[t]
\centering
\footnotesize
\setlength{\tabcolsep}{4pt}
\renewcommand{\arraystretch}{1.08}
\begin{tabularx}{0.96\textwidth}{@{}
>{\raggedright\arraybackslash}p{3.3cm}
>{\raggedright\arraybackslash}p{4.8cm}
>{\raggedright\arraybackslash}p{3.5cm}
>{\raggedright\arraybackslash}p{2cm}
>{\raggedright\arraybackslash}X
@{}}
\toprule
\textbf{Failure mode} &
\textbf{Root cause} &
\textbf{Best signal} &
\textbf{Evidence class} &
\textbf{Typical fix} \\
\midrule
Interface / shape mismatch
& Wrong tensor shape or missing decomposition output
& Local validation
& Neither (local)
& CodeFix \\

\addlinespace[0.25em]
Wrong incentive response
& Sign error, missing coupling term, or mis-specified objective
& Private-objective sign (C1) + price-response (C2)
& Behavioral
& Reformulate \\

\addlinespace[0.25em]
Cost attached to fixed data
& Cost multiplies a parameter rather than a decision variable
& Pattern check (C6) + price-response (C2)
& Both
& Reformulate \\

\addlinespace[0.25em]
Missing cost term
& Cost key exists in data but unused in code
& Missing-cost (C11)
& Static
& CodeFix or Reformulate \\

\addlinespace[0.25em]
Constraint scope bug
& Missing index, incorrect aggregation, or omitted constraint
& Convergence (C3) + decision-variable audit (C7)
& Both
& Reformulate \\

\addlinespace[0.25em]
Degenerate proposals
& All-zero plan, dominated objective, or bad bounds
& Degeneracy (C4)
& Behavioral
& Reformulate \\
\bottomrule
\end{tabularx}
\caption{Failure-mode taxonomy for generated optimization agents. Check numbers refer to the evidence checks in Table~\ref{tab:checks_list}. The Evidence class column indicates whether the failure is primarily caught by behavioral (in-loop) checks, static (code-level) checks, or both, corresponding to the ablation in Section~\ref{subsec:evidence_ablation}.}
\label{tab:taxonomy}
\end{table*}

The three scenario-level categories discussed in Section~\ref{subsec:failure_analysis} arise from different bug-level mechanisms in Table~\ref{tab:taxonomy}. Novel cost structures most often lead to missing or mis-specified objective terms. Ambiguous specifications can produce wrong incentive response or constraint-scope bugs because multiple formalizations are plausible from the text. Complex constraint composition primarily appears as scope, indexing, or aggregation errors. Together, these cases explain why the remaining failures are mostly structural rather than syntactic.

\section{Worked Scenario Examples}
\label{app:worked_examples}

We provide three representative scenarios spanning the difficulty range: an easy training scenario, a medium test scenario that illustrates a failure mode caught by coordination-in-the-loop verification, and a hard test scenario with novel cost structure.

\subsection{Example 1: Training Scenario (Easy)}
\label{app:example_easy}

\paragraph{Tags.} \texttt{QP, order\_aggregation, warehouse\_capacity, transport\_cost}

\paragraph{Natural-language specification.}
``A vendor operates a single warehouse and ships 2 ASINs to 2 Retail inbound nodes over 3 weeks. The warehouse has a maximum outbound shipping capacity of 150 units per ASIN per week. Transportation costs vary by destination: to node 1, \$2 per unit; to node 2, \$3 per unit. The vendor wants to minimize total transportation cost while respecting the warehouse outbound capacity.''

\paragraph{Data schema.}
\begin{itemize}
    \item \texttt{transport\_cost}: shape $(J,)$
    \item \texttt{warehouse\_capacity}: float
    \item Problem size: $A=2$, $J=2$, $T=3$
\end{itemize}

\paragraph{Intended formulation.}
The agent defines shipment variables $s_{i,w,j,t} \ge 0$ (or equivalently $s_{i,j,t}$ in the single-warehouse case), enforces order aggregation
\begin{equation}
\mathrm{po}_{i,j,t} = \sum_{w} s_{i,w,j,t},
\end{equation}
and imposes the outbound-capacity constraint
\begin{equation}
\sum_{j} s_{i,w,j,t} \le \mathrm{warehouse\_capacity}
\qquad
\forall i,t.
\end{equation}
The objective is
\begin{equation}
\min \sum_{i,j,t} \mathrm{transport\_cost}_{j} \cdot s_{i,j,t}.
\end{equation}

\paragraph{Why this scenario is easy.}
This scenario uses only familiar linear costs and constraints from the OR block library. The main modeling pitfall is attaching transportation cost to a fixed input parameter rather than to the shipment decision variable --- a bug that Check~6 (cost$\times$data anti-pattern) is designed to catch.

\subsection{Example 2: Test Scenario (Medium) --- Failure Caught Only by Behavioral Evidence}
\label{app:example_medium}

\paragraph{Tags.} \texttt{QP, order\_aggregation, inventory\_balance, shared\_capacity, ambiguous\_nl}

\paragraph{Natural-language specification.}
``Our operation is straightforward, but the capacity situation is a bit nuanced. We have two warehouses, three products, two Retail docks, and three weeks. Making products A, B, and C costs \$0.15, \$0.25, and \$0.20 per unit, respectively. Shipping costs \$0.10 per unit on every lane. Holding inventory costs \$0.20 per unit per week. Each warehouse starts with 40 units of every product and receives 25 units of each product each week. The catch is that products A and B share a production line: together, they cannot exceed 100 units per warehouse per week. Product C has its own line with capacity 80 units per warehouse per week. In addition, each warehouse can ship at most 150 units per week in total across all products and destinations. Keep costs down.''

\paragraph{Data schema.}
\begin{itemize}
    \item \texttt{production\_cost}: shape $(A,)$
    \item \texttt{transport\_cost\_per\_unit}: float
    \item \texttt{holding\_cost}: float
    \item \texttt{initial\_inventory}: shape $(A, W)$
    \item \texttt{procurement}: shape $(A, W, T)$
    \item \texttt{shared\_line\_capacity}: float
    \item \texttt{dedicated\_line\_capacity}: float
    \item \texttt{outbound\_capacity}: float
    \item Problem size: $A=3$, $W=2$, $J=2$, $T=3$
\end{itemize}

\paragraph{Common LLM failure.}
The specification says that each warehouse receives 25 units of each product every week. This is a fixed inflow. However, because the text also says that ``making products'' has a per-unit cost, the LLM may introduce an extra controllable production variable $\mathrm{prod}[w,a,t]$. The generated inventory balance then becomes $\mathrm{inv} = \mathrm{prev\_inv} + \mathrm{prod} + \mathrm{procurement} - \mathrm{shipments}$. This is wrong because the model now has two sources of supply: the fixed procurement inflow and the new production decision. Since production is costly, the optimizer can set $\mathrm{prod}=0$ and still use the fixed procurement inflow to satisfy orders. The model compiles and returns feasible purchase-order quantities, but the economic meaning is wrong: fulfilled units are no longer charged the intended production cost.

\paragraph{Why static checks miss it.}
From the code alone, the formulation looks reasonable. It has an inventory balance, it has production-capacity constraints, and it has a production-cost term in the objective. The variable $\mathrm{prod}$ also appears plausible because the natural-language specification mentions both ``making products'' and shared production capacity. A static checker can flag the variable for review, but the diagnosis LLM may still accept it as intended. The key error is not a missing constraint or a syntax failure. The key error is that $\mathrm{prod}$ interacts with the fixed procurement inflow to create a free source of supply. That interaction is hard to detect by inspecting individual code fragments.

\paragraph{How behavioral checks catch it.}
Check~8 (marginal-cost consistency) solves the agent at two consensus levels ($z{=}30$ and $z{=}80$) and computes the marginal change in utility per unit of PO. For a cost-minimizing vendor, asking it to ship more units should make reported private cost increase, because additional PO volume should create additional production, shipping, holding, or capacity cost. In other words, if the vendor is asked to do more work, the model should not report that the vendor is better off. Instead, the generated model has a \emph{positive} marginal: increasing PO increases reported utility. This happens because the optimizer sets $\mathrm{prod}=0$ to avoid production cost and absorbs additional volume through the fixed procurement inflow, so the production-cost term contributes nothing to fulfilled units. This wrong-sign marginal is a clear behavioral signal that the private cost function is structurally incomplete. The behavioral-only diagnosis LLM therefore triggers \textsc{Reformulate}, producing a corrected formulation that removes the spurious production variable and applies production cost to the fulfilled flow.

\paragraph{Result.}
The wrong-sign behavioral signal triggers \textsc{Reformulate}. The corrected formulation removes the spurious production variable and attaches the per-unit production cost to the fulfilled flow. In our ablation, behavioral-only verification fixes this scenario and reaches a 0\% objective gap, while static-only verification fails with an 89.9\% gap. This example shows why behavioral verification is useful: the incorrect model is executable, feasible, and locally plausible, but its response to a higher purchase-order target reveals that its economics are wrong.

\subsection{Example 3: Test Scenario (Hard)}
\label{app:example_hard}

\paragraph{Tags.} \texttt{QP, order\_aggregation, inventory\_balance, congestion\_cost, ambiguous\_nl}

\paragraph{Natural-language specification.}
``We run a couple of warehouses and move two product lines to Retail through two receiving docks over three weeks. Our stuff costs about \$0.30 per unit to make. We keep track of what's in the warehouse---at the start of the planning window we've got 50 of everything everywhere. Each warehouse gets 30 units of each product per week from our suppliers. Storage runs us \$0.40 per unit sitting around each week. Shipping is \$1.50 per unit base, but when a lane gets busy the cost goes up---there's a congestion factor of 0.02 per unit squared on top of the base rate.''

\paragraph{Data schema.}
\begin{itemize}
    \item \texttt{production\_cost}: float
    \item \texttt{initial\_inventory}: shape $(A, W)$
    \item \texttt{procurement}: shape $(A, W, T)$
    \item \texttt{holding\_cost}: float
    \item \texttt{base\_transport\_cost}: float
    \item \texttt{congestion\_factor}: float
    \item Problem size: $A=2$, $W=2$, $J=2$, $T=3$
\end{itemize}

\paragraph{Intended formulation.}
The correct model requires shipment variables, inventory-balance constraints, and order aggregation:
\begin{equation}
\mathrm{po}_{i,j,t} = \sum_{w} s_{i,w,j,t}.
\end{equation}
Inventory evolves according to standard balance equations using initial inventory, weekly procurement, outbound shipments, and ending inventory. The transportation term includes both a linear base cost and a quadratic congestion term applied to the shipment decision:
\begin{equation}
\sum_{i,w,j,t}
\left(
\mathrm{base\_transport\_cost} \cdot s_{i,w,j,t}
+
\mathrm{congestion\_factor} \cdot s_{i,w,j,t}^{2}
\right)
\end{equation}
The full objective also includes production and holding costs.

\paragraph{Common LLM failure.}
The LLM often treats \texttt{congestion\_factor} as a fixed surcharge or multiplies it by \texttt{base\_transport\_cost} (a data$\times$data pattern) rather than applying it as a quadratic term on the shipment decision variable. The informal phrasing (``gets busy,'' ``on top of the base rate'') makes the correct interpretation ambiguous.

\paragraph{Why this scenario is hard.}
This scenario combines informal natural-language phrasing with a nonlinear transport-cost structure. Both the \texttt{ambiguous\_nl} and \texttt{congestion\_cost} tags are absent from the training split, so episodic memory provides no direct precedent. The agent must infer that the congestion factor defines a per-lane quadratic cost on the shipment decision variable. Even with full evidence, this scenario fails in approximately 30\% of runs, placing it among the residual hard cases discussed in Section~\ref{subsec:failure_analysis}.

\end{document}